\documentclass{ifacconf}

\usepackage{graphicx}      
\usepackage{natbib}        
\usepackage{amssymb}
\usepackage{amsmath}
\usepackage{mathtools}
\usepackage{bbm}
\usepackage{enumerate}
\usepackage{xfrac}

\newcommand{\D}{\mathcal{D}}
\newcommand{\F}{\mathcal{F}}

\renewcommand{\L}{\mathcal{L}}
\newcommand{\la}{\lambda}

\newcommand{\R}{\mathbb{R}}
\newcommand{\X}{\mathcal{X}}
\renewcommand{\d}{\mathrm{d}}

\begin{document}
\begin{frontmatter}

\title{Port-Hamiltonian systems with energy and power ports\thanksref{footnoteinfo}} 

\thanks[footnoteinfo]{Kaja Krhač acknowledges  funding from the European Union Horizon Europe MSCA Grant No. 101073558 (ModConFlex) and the IMPACTS project ANR-21-CE48-0018}

\author[First]{Kaja Krhač} 
\author[Second]{Bernhard Maschke} 
\author[Third]{Arjan van der Schaft}

\address[First]{Functional Analysis Group, University of Wuppertal, Germany (e-mail: krhac@uni-wuppertal.de).}
\address[Second]{LAGEPP, Claude Bernard University Lyon 1, France \\(e-mail: bernhard.maschke@univ-lyon1.fr)}
\address[Third]{Bernoulli Institute for Mathematics, Computer Science and AI, 
   University of Groningen, the Netherlands, \\(e-mail: a.j.van.der.schaft@rug.nl)}

\begin{abstract}
We extend the port-Hamiltonian framework defined with respect to a Lagrangian submanifold and a Dirac structure by augmenting the Lagrangian submanifold with the space of external variables. The new pair of conjugated variables is called energy port.
We show that in the most general case, the extension describes constrained Hamiltonian systems whose Hamiltonian function depends on inputs.

\end{abstract}
\begin{keyword}
Port-Hamiltonian systems, constrained Hamiltonian systems, input-output systems, Morse families, Lagrangian submanifolds
\end{keyword}

\end{frontmatter}

\section{Introduction}
In this paper we extend the finite-dimensional port-Hamiltonian framework defined on a Lagrangian submanifold of the cotangent bundle of a state space~(\cite{DAE},~\cite{vietnam}) by introducing a pair of conjugated variables on the Lagrangian submanifold that we call energy ports. Unlike ports in the Dirac structure, referred to as power ports, energy ports are \emph{not} dual variables. They arise naturally in infinite-dimensional systems~(\cite{maschke2023linear}) as boundary ports, whereas in this paper we introduce them by augmenting the state space with the space of external variables which in turn augments the cotangent bundle. Using the description of Lagrangian submanifolds given by Morse families, we show that the new definition allows us to describe constrained Hamiltonian systems~(\cite{dirac2001lectures}) whose Hamiltonian depends on external variables called inputs and whose outputs are quantities that are conjugated with respect to the given inputs~(\cite{asFeedback}).

Port-Hamiltonian theory generalizes the Hamiltonian description for closed systems to open systems by introducing so-called port variables or ports. The energy of a port-Hamiltonian system is given by the Hamiltonian function, a smooth function on a manifold called the state space of the system. Ports, on the other hand, capture the rate of change of energy of an open system in time, since they are chosen as dual variables whose duality pairing gives the power. We will refer to such ports as \emph{power} ports in order to distinguish them from newly defined \emph{energy} ports in Sec.~\ref{sec:results}.  Energy ports are not dual variables, but they are external variables that appear in the duality pairing of power ports. 
Moreover, ports can be used to interconnect open systems. Energy  then flows between different components of the total system in such a way that the total power is conserved. 
The mathematical structure which captures this power-conserving interconnection between ports and in turn determines the dynamics of the system is a so-called Dirac structure. Further generalization of port-Hamiltonian systems to Lagrangian submanifolds concerns the role played by the Hamiltonian function in the dynamics of the system. The generalization
identifies the Lagrangian submanifold as the geometric structure which together with the Dirac structure governs the dynamics of the system. As we recall in Sec.~\ref{sec:revision}, 
Lagrangian submanifolds of the cotangent bundle can be described by multi-valued function on the base manifold, called Morse family. As a result of this correspondence, Lagrangian submanifold is the structure which allows us to lift the characterization of energy from the state space to its cotangent bundle. The implication is two-fold. On the one hand, the multi-valued nature allows us to describe constrained Hamiltonian systems in Sec.~\ref{sec:results}. These systems admit a parameterized \emph{family} of Hamiltonian functions where each member is a viable Hamiltonian function from the physical standpoint. On the other hand, the Lagrangian submanifold can be straightforwardly extended with a space of external variables by extending the state space. These additional variables defined in Sec.~\ref{sec:results} are a conjugated pair of variables called energy ports.
The extension allows us direct access to external variables appearing in the Hamiltonian function. 





\emph{Notation and conventions:} Throughout the paper we use Einstein summation convention: indices that repeat are summed over. We use the same notation for coordinates and points of a manifold. Elements of the dual vector space of some vector space are called co-vectors. 
If not specified, a Lagrangian submanifold is a Lagrangian submanifold of the cotangent bundle of a smooth manifold. The symplectic structure on the cotangent bundle is the canonical symplectic 2-form.

\section{Port-Hamiltonian systems on Lagrangian submanifolds}\label{sec:revision}
In this section we briefly recall definitions and results regarding Lagrangian submanifolds of a cotangent bundle equipped with the canonical symplectic form that are relevant for understanding port-Hamiltonian systems defined on the Lagrangian submanifold. 

In~(\cite{vietnam}), a port-Hamiltonian system on a Lagrangian submanifold is defined as follows. 
\begin{defn}\label{def:phs}
    A port-Hamiltonian system defined on a Lagrangian submanifold is a quadruple
    $$
    (\X, \F^p, \L, \D)
    $$
    consisting of a smooth $n$-dimensional manifold $\X$, a smooth vector bundle $\F^p$ over $\X$ of rank $d$ and its dual bundle $\F^{p*}$, a Lagrangian submanifold $\L \subset T^*\X$ of the cotangent bundle $(T^*\X, \omega)$ equipped with the canonical symplectic 2-form $\omega$ and a Dirac structure $\D$ of the Whitney sum of ${T\X\oplus \F^p\oplus T^*\X\oplus\F^{p*}}$  with respect to duality pairing 
    $$
    \langle \cdot | \cdot \rangle + \langle \cdot | \cdot \rangle_p
    $$
    where $\langle \cdot | \cdot \rangle$ denotes the duality pairing between  $T_x\X$ and $T^*_x\X$ and $\langle \cdot | \cdot \rangle_p$ denotes duality pairing between $\F^p_x$ and $\F^{p*}_x$, for any $x\in\X$. The time evolution of a physical system modelled by $(\X, \F^p, \L, \D)$ is given by a differentiable curve $x:I\subseteq\R \to \X$ such that
\begin{gather}
        (x, e)(t) \in \L   \label{eq:L} \;, \quad 
        (\dot x, f^p, e, e^p)(t) \in \D
\end{gather}
for every $t \in I$, where $\dot x (t)$ denotes the tangent at $x(t)$.
\end{defn}
The manifold $\X$ is called the state space of a system. 
Pairs of variables $(f^p, e^p)\in \F_x^p \oplus \F_x^{p*}$ are dual variables called open power ports. They allow one to connect another open system with the same port space. The duality pairing of open port variables along the trajectory of the system encodes energy flow, or power, between the system and the environment. From the power-conserving property of the Dirac structure,
$$
\langle e(t) \mid \dot x(t) \rangle + \langle e^p(t) \mid f^p(t) \rangle_p = 0 \;
$$
for any $(\dot x, f^p, e, e^p)(t) \in \D$.
On the other hand, a closed system has no open power ports. The Dirac structure $\D^c$ of a closed system is defined on ${T\X\oplus T^*\X}$ with respect to duality pairing $\langle \cdot | \cdot \rangle$. There is no exchange of energy with the environment and the total energy is conserved, 
$$
\langle e(t) \mid \dot x(t) \rangle  = 0 \;
$$
for any $(\dot x, e)(t) \in \D^c$.
In order to determine time evolution of the system, one 
must first be able to find a co-vector at every point of the curve, such that all pairs constitute the Lagrangian submanifold $\L$.  The tangent vector at every point of the curve, together with the co-vector dictated by the Lagrangian submanifold and the open power ports must belong to the given Dirac structure. 

We now show how a physical systems whose energy is described by the Hamiltonian function fits into the Def.~\ref{def:phs}. The theorem proved by Maslov~(\cite{maslov1972théorie}) and refined by Hörmander~(\cite{Hormander1971-je}) shows that any Lagrangian submanifold of the cotangent bundle equipped with the canonical symplectic 2-form, in particular, can be described in a neighbourhood of each of its points by a parameterized family of functions called Morse family.
The family can consist of a single function, which, in the context of port-Hamiltonian framework, describes systems with the Hamiltonian function. The result also prepares the ground for energy ports introduced in the next section.

The statement of the theorem is in so-called canonical coordinates.
Recall that $\omega= \d \theta$ where $\theta$ is the canonical {1-form}. 
Let ${x=(x_1,\dots, x_n)\in \R^n}$ denote local coordinates on $\X$ inducing coordinates ${e=(e_1, \dots, \e_n) \in \R^n}$ on the fiber of $T^*\X$. Then $(x,e)\in \R^{2n}$ are local coordinates on $T^*\X$ in which $\theta$ and $\omega$  have the form
$$
\theta = \sum_{i=1}^n e_i\d x^i \;, \quad \omega = \d\theta=  \sum_{i=1}^n \d e_i \wedge \d x_i \;.
$$
We call such coordinates canonical coordinates. Coordinates $x$ on $\X$ are called energy variables, as the Hamiltonian function is a function on $\X$. Coordinates $e$ on the fiber of $T^*\X$ are then called co-energy variables, since they are components of a co-vector on the fiber. The pair $(x,e)$ are called \emph{conjugated variables}.

\begin{thm}\label{th:MH} (Maslov-Hörmander)
  Let $K:\X \times \R^k \to \R$ be a family of smooth functions on $\X$ parameterized by $k\geq0$ parameters.  Let $x=(x_1,\dots, x_n)\in \R^n$ denote local coordinates on $\X$, $\la = (\la_1, \dots, \la_k) \in \R^k$ parameters and $K(x, \la)$ the real value of $K$ in the local coordinates. We define critical set $\Sigma_K$ of $K$ as
  $$
    \Sigma_K \coloneqq \left\{ 
    (x, \la) \in \R^n \times \R^k 
    \mid  
    \frac{\partial K}{\partial \la^i}(x, \la) = 0 \;, \; i=1,\dots, k
    \right\}\;.
  $$
If  the  condition
\begin{equation}\label{eq:rank}
\mathrm{rank}
\begin{pmatrix}
    \frac{\partial^2K}{\partial \la^i \partial x^j}(x,\la) & 
    \frac{\partial^2K}{\partial \la^i \partial \la^j}(x,\la)
\end{pmatrix}\bigg|_{\Sigma_K} = k
\end{equation} 
holds, i.e., if the $k \times (n+k)$ matrix has maximal rank at every point $(x,\la)\in \Sigma_K$, then the subset $\L_K \subset T^*\X$, given in canonical coordinates as
\begin{equation}\label{eq:LMorse}
\L_K \coloneqq \left\{
(x,e) \in \R^{2n} \mid 
\exists \la \in \R^k: e = \frac{\partial K}{ \partial x}(x,\la)\bigg|_{\Sigma_K}
\right\}
\end{equation}
is a Lagrangian submanifold of $(\R^{2n}, \sum_{i=1}^n \d e_i \wedge \d x_i)$.
Vice versa, if $\L$ is a Lagrangian submanifold of $(T^*\X, \omega)$ equipped with the canonical symplectic 2-form $\omega$, then, in the neighbourhood of every point in $\L$, for some $k\geq 0$, $\L$ can be described by Eq.~\ref{eq:LMorse} for some $K$ that satisfies Eq.~\ref{eq:rank}.
\end{thm}
\begin{defn}\label{def:gfMf}
    A  family of functions $K:\X\times \R^k\to \R$ on $\X$ parameterized by $k>0$ parameters satisfying Eq.~\ref{eq:rank} is called \emph{Morse family}. If $k=0$, then $K:\X\to\R$ is called a \emph{generating function} of the Lagrangian submanifold. We say that $\L_K$ is generated by a Morse family $K$, or a generating function if $k=0$.
\end{defn}

The proof of the theorem can be found in~(\cite{weinstein}, \cite{cardin2014elementary}). Here we only prove sufficiency for $k=0$. This case describes Lagrangian submanifolds generated by a single function. 
\begin{prop}\label{prop:Ltrans} Consider the set up of Th.~\ref{th:MH} with $k=0$. 
    Let $K:\X \to \R$ be a smooth function. Then $\L_K \subset T^*\X$, given in canonical coordinates by Eq.~\ref{eq:LMorse}:
        $$ 
    \L_K  \coloneqq  \left\{ 
    (x, e) \in \R^{2n} \mid e_i = \frac{\partial K}{\partial x^i}(x) \;, \; i=1, \dots, n
    \right\}
    $$
    is a Lagrangian submanifold of $(\R^{2n}, \sum_{i=1}^n \d e_i \wedge \d x_i)$.
\end{prop}
\begin{pf}
    Note that $\L_K$ is a coordinate representation of the image $\d K_\X$ of the 1-form $\d K: \X \to T^*\X$.  A point $(x,e) \in \d K_\X$ is a co-vector $e \in T_x^*\X$ with the base point $x \in \X$ such that $e=\d K(x)$.
    The components $(e_1, \dots, e_n)$ of $e$ in canonical coordinates are $e_i = \partial_i K$ for every $i=1,\dots, n$.
    To show $\d K_\X$ is a Lagrangian submanifold, note that a 1-form is a canonical embedding by definition,
    $$
    \dim (\d K_\X) = \dim \X =  \frac{1}{2}\dim T^*\X 
    $$
    The submanifold $\d K_\X$ is therefore Lagrangian with respect to ${\omega=\d\theta}$,
    $$
    (\d K)^*\d \theta = \d((\d K)^{*} \theta) = \d(\d K) = 0 
    $$
    The first equality follows because pull-back commutes with exterior derivative. The second equality is a consequence of tautological property of the canonical 1-form, ${(\d K)^*\theta=\d K}$ and the last equality follows because $\d^2=0$. 
\end{pf}
Geometrically, the image of a differential of a function generates a Lagrangian submanifold $\L$ that is at each of its points transversal to the fiber of $T^*\X$ that passes through that point and $\L \to \X$ is a diffeomorphism. Since the opposite direction also holds~(\cite{libermann_marle_1989}), one can equivalently characterize Lagrangian submanifolds generated by a function as transversal. Non-transversal Lagrangian submanifolds are generated by Morse families.

In view of Th.~\ref{th:MH}, a Lagrangian submanifold is  a structure which describes relations among energy and co-energy variables. We call these relations \emph{constitutive relations}. 
Physical systems whose energy is given by the Hamiltonian function $H:\X\to\R$   
can be understood as a class of port-Hamiltonian systems with transversal Lagrangian submanifold $\L_H = \d H_\X$.
The  condition~(\ref{eq:L}) in Def.~\ref{def:phs} reduces to $e(t)=\d H(x(t))$ for any $t \in I$
and one can recover the familiar definition of a port-Hamiltonian system as a quadruple $(\X, \F^p, H, \D)$ with dynamics of the system given by the differential equation~(\cite{intro})
$$
(\dot x(t), f^p(t), \d H(x(t)), e^p(t)) \in \D
$$
at every $t\in I$. The generating function $H$ plays only an auxiliary role in determining the dynamics of the system, since only the differential of $H$ appears in the differential equation. It is the Lagrangian submanifold $\L_H$ which $H$ generates that, together with the Dirac structure, determines the dynamics. 
As constitutive relations determined by transversal $\L_H$ are also a mapping, one can understand the generalization of port-Hamiltonian to Lagrangian submanifolds as a way to describe systems with relations among energy and co-energy variables that are not given by a map. It is this multi-valued nature of a Lagrangian submanifold what makes it possible to describe the total Hamiltonian of constrained Hamiltonian systems (see Sec.~\ref{sec:results}).
On the other hand, the generator of the Lagrangian submanifold, whether it is a function or a Morse family, is related to power balance. Along the trajectory of a port-Hamiltonian system $(\X, \F^p, \L_H, \D )$,
\begin{equation}\label{eq:pb}
    \frac{\d H(x(t), \la(t))}{\d t} + \langle e^p | f^p \rangle_p = 0 
\end{equation}
When we introduce energy ports in Sec.~\ref{sec:results} the power balance will not be given with respect to the generator of the Lagrangian submanifold.

\section{Hamiltonian systems with constraints and energy ports}\label{sec:results}

We further exploit the structure of Lagrangian submanifolds by introducing variables, called energy ports, which account for the change in energy of an open systems caused by a change in constitutive relations. 
A port-Hamiltonian system defined with respect to such Lagrangian submanifold will, in general, have ports of two different types: one associated with the constitutive relations (energy ports) and one associated with the interconnection structure of the physical system (power ports). Unlike power ports, energy ports are not dual variables but they play a role in power balance, which justifies their name. 

Energy ports are defined with respect to a specific class of Morse families. 
\begin{defn}
    Let $K:\X \times \R^m \times \R^k \to \R$ be a family of functions on $\X \times \R^m$ parameterized by $k \geq 0$ parameters, such that $K$ is a Morse family if $k>0$. Let $\L_K$ be the Lagrangian submanifold of $(T^*(\X \times \R^m), \d \theta)$ generated by $K$,
\begin{equation}\label{eq:LEP}
\begin{gathered}
\L_K \coloneqq \bigg\{
(x, \nu^p, e, \varepsilon^p) \in \R^{2(n+m)} \mid \exists \lambda \in \R^k :  \\
e = \frac{\partial K}{\partial x}(x, \nu^p, \la)\bigg|_{\Sigma_K} \;, \quad 
\varepsilon^p = \frac{\partial K}{\partial \nu^p}(x, \nu^p, \la)\bigg|_{\Sigma_K}  \bigg\}
\end{gathered}
\end{equation}
where  $(x, \nu^p) \in \R^{n+m}$ denotes local coordinates on ${\X\times\R^m}$. The induced local coordinates on the fiber of ${T^*(\X \times \R^m)}$ are $(e, \varepsilon^p) \in \R^{n+m}$, so that the canonical symplectic 2-form ${\d\theta=\sum_{i=1}^{n} \d e_i\wedge\d x_i + \sum_{i=1}^{m} \d \varepsilon^p_i\wedge\d \nu^p_i }$ and  $\Sigma_K$ is the critical set of $K$, 
$$
\Sigma_K \coloneqq \left\{ 
(x, \nu^p, \lambda) \in \R^{n}\times\R^m\times\R^k \mid \frac{\partial K}{\partial \lambda} (x, \nu^p, \lambda) = 0 
\right\} \;.
$$
If it holds that
\begin{equation}\label{eq:rankEP}
\mathrm{rank}
\begin{pmatrix}
    \frac{\partial^2K}{\partial \nu^{pi} \partial x^{pj}}(x, \nu^p, \la) & 
    \frac{\partial^2K}{\partial \nu^{pi} \partial \nu^{pj}}(x, \nu^p, \la)
\end{pmatrix}\bigg|_{\Sigma_K} = m 
\end{equation} 
that is, if $m\times(n+m)$ matrix has maximal rank at every point of the critical set $\Sigma_K$, then for any 
$(x, \nu^p , e, \varepsilon^p) \in \L_K$ a pair of variables $(\nu^p, \varepsilon^p) \in T^*\R^m$ is an open \emph{energy port} and $\L_K$ is called \emph{port-Lagrangian submanifold}.
\end{defn}
We point out that derivatives in the rank condition~(\ref{eq:rankEP}) are only with respect to coordinates on the base space $\X\times \R^m$, whereas~(\ref{eq:rank}) concerns parameters in $\R^k$. Thus, a Morse family on 
 $\X\times \R^m$ that satisfies~(\ref{eq:rankEP})
 generates the Lagrangian submanifold~(\ref{eq:LEP}), according to Th.~\ref{th:MH}. 
Moreover, because of the local character of Th.~\ref{th:MH}, the space of parameters $\R^k$ and the base space of energy ports $\R^m$ can both be replaced by more complicated smooth manifolds. 
\begin{defn}
    A Morse family $K:\X\times\R^m\times\R^k\to\R$ on ${\X\times\R^m}$ parameterized by $k>0$ parameters that satisfies Eq.~(\ref{eq:rankEP}) is called \emph{restricted Morse family}. If $k=0$ and $K:\X\times\R^m\to\R$ satisfies Eq.~(\ref{eq:rankEP}) we call $K$ a \emph{restricted generating function}. 
\end{defn}
We define a port-Hamiltonian system with respect to a port-Lagrangian submanifold as follows. 
\begin{defn}
  Let $\X$ be a smooth $n$-dimensional manifold, $\F^p$ a smooth vector bundle over $\X$ of rank $d$ and $\F^{p*}$ its dual bundle. A port-Hamiltonian system defined on a port-Lagrangian submanifold is a quintuple
  $$
  (\X, \F^p, \R^m, \L^p, \D)
  $$
  where $\L^p\subset T^*(\X\times\R^m)$ is a port-Lagrangian submanifold of the cotangent bundle $(T^*(\X\times\R^m), \omega)$ equipped with the canonical symplectic 2-form $\omega$ and $\D$ is a Dirac structure ${\D\subset T\X\oplus\F^p\oplus T^*\X\oplus\F^{p*}}$ with respect to duality pairing
  $$
    \langle \cdot | \cdot \rangle + \langle \cdot | \cdot \rangle_p
    $$
    where $\langle \cdot | \cdot \rangle$ denotes the duality pairing between  $T_x\X$ and $T^*_x\X$ and $\langle \cdot | \cdot \rangle_p$ denotes duality pairing between $\F^p_x$ and $\F^{p*}_x$, for any $x\in\X$. The time evolution of a physical system described by a port-Hamiltonian system $  (\X, \F^p, \R^m, \L^p, \D)$ is given by a curve $x:I\subseteq \R \to \X$ satisfying
\begin{gather}
        (x, \nu^p, e, \varepsilon^p)(t) \in \L^p   \label{eq:Ldep} \\
        (\dot x, f^p, e, e^p)(t) \in \D \label{eq:Ddep}
\end{gather}
for any ${t\in I}$, where $\dot x(t)$ denotes the tangent at $x(t)$.
\end{defn}
Note that the Dirac structure is a vector subbundle over $\X$, whereas the Lagrangian submanifold is defined over ${\X\times\R^m}$. The energy ports will, however, appear in the Dirac structure, and consequently in the power balance. The pair of conditions is coupled, so
$e$ will be given in terms of relations among $x, \nu^p, \la$ and $\varepsilon^p$  as dictated by the port-Lagrangian submanifold $\L^p$. 
$$
\langle e(t) | \dot x(t) \rangle + \langle e^p(t) | f^p(t) \rangle_p = 0 \;
$$
In local coordinates in a neighbourhood of $x(t)$ for any $t\in I$, the power balance can be written as
$$
\frac{\partial H_T}{\partial x^i}(x(t), \nu^p(t), \lambda(t))\dot x^i(t) + e^{p}_jf^{pj} = 0 
$$
where $H_T:\X\times\R^m\times\R^k\to\R$ is the restricted Morse family on $\X\times \R^m$ that generates $\L^p$ (see Eq.~\ref{eq:Ldep}) and $\la \in \R^k$.
This power balance is not with respect to $H_T$, since
\begin{gather*}
    \frac{\d H_T}{\d t} = \frac{\partial H_T}{\partial x^i}\dot x^i 
    +\frac{\partial H_T}{\partial \nu^{pi}}\dot \nu^{pi} + 0
\end{gather*}
The expression is evaluated at $(x(t), \nu^p(t), \la(t))$. We henceforth suppress the dependencies to make the notation less cluttered.  
Since $\nu^p$ are the independent variables whose dynamics we do not consider, we use the chain rule to eliminate $\dot \nu^p$, which leads to
$$
 \frac{\partial H_T}{\partial x^i}\dot x^i  = \frac{\d }{\d t}\left( 
H_T - \frac{\partial H_T}{\partial \nu^{pi}}\nu^{pi}
\right) + \left( \frac{\d}{\d t} \frac{\partial H_T}{\partial \nu^{pi}} \right)\nu^{pi}
$$
This expression only makes sense in coordinates. The total power balance can then be written as 
$$
\frac{\d \widetilde H}{\d t} + \frac{\d \varepsilon^p_i}{\d t} \nu^{pi} + e^{p}_jf^{pj} = 0
$$
where we introduced a function $\widetilde H: \X\times\R^m\times\R^k\to\R$, 
$$
\widetilde H (x, \nu^p, \la) \coloneqq H_T(x, \nu^p, \la) - \frac{\partial H_T}{\partial \nu^{pi}}(x, \nu, \la)\nu^{pi}
$$
As expected, the power balance has an additional term compared to~(\ref{eq:pb}). The second term describes energy change due to a change in constitutive relations and  the dual variable to $\nu^p$ is $\dot\varepsilon^p$~(\cite{MASCHKE1992}).

We now consider special cases of port-Hamiltonian systems on port-Lagrangian submanifolds that are generated by a restricted Morse family ${H_T:\X\times\R^m\times\R^k\to\R}$ on ${\X\times\R^m}$. 
We show that ${m=0},k>0$ describes a constrained Hamiltonian system. Such port-Hamiltonian systems do not have energy ports and are described in Def.~\ref{def:phs}. 
On the other hand, $m>0, k=0$ describes input-output Hamiltonian systems without constraints. Its Lagrangian submanifold is generated by a restricted generating function on $\X\times\R^m$. We will more closely study an input-output Hamiltonian system that is linear in inputs. 
If $m=k=0$, then $H_T$ is a generating function on $\X$. 
Finally, we put everything together and give an example of a constrained input-output Hamiltonian system that is non-linear in inputs.

\subsection{Input - output Hamiltonian systems}
Input-output Hamiltonian systems are open physical systems whose constitutive relations depend on variables whose dynamics we do not consider. Such variables are then viewed as external variables and called inputs. 

We study more closely an input-output Hamiltonian system that is linear in inputs. Let $Q$ be a smooth $d$-dimensional manifold. Consider $H_T:T^*Q\times\R^m\to\R$ defined as
\begin{gather*}
H_T(x, u) = H(x)+u^iG_i(x) \\
y_i = G_i(x) \;, \quad i=1,\dots, m
\end{gather*}
where $G_i:T^*Q\to\R$ for every $i=1,\dots, m$ and ${u=(u_1,\dots, u_m)\in\R^m}$ are inputs.  We assume that $H_T$ is a restricted generating function on $T^*Q\times \R^m$. That is, at every point $(x,u)\in \R^{2d}\times\R^m$
$$
\mathrm{rank}
\begin{pmatrix}
    \frac{\partial G_1}{\partial x}(x) \\
    \vdots \\
    \frac{\partial G_m}{\partial x}(x) \\
\end{pmatrix} =m
$$
Thus, condition~(\ref{eq:rankEP}) in this case ensures that there is no redundancy in inputs. The port-Lagrangian submanifold $\L^p$ of $T^*T^*Q$ it generates is 
\begin{gather*}
\L^p = \bigg\{ 
(x,u, e, y) \in R^{2(2d+m)} \mid  \\
e = \frac{\partial H_T}{\partial x} (x,u) =  \frac{\partial H}{\partial x}(x) + u^i \frac{\partial G^i}{\partial x}(x) \\
y_i= \frac{\partial H_T}{\partial u} (x,u) = G_i(x)\;, \quad i=1, \dots, m
\bigg\}
\end{gather*}
Such input-output system is described by a port-Hamiltonian system $(T^*Q, \R^m, \L^p, \D^c)$ where $\D^c\subset TT^*Q\oplus T^*T^*Q$ is a closed Dirac structure defined as the graph of the Poisson vector field on $(T^*T^*Q, \omega)$ where $\omega$ is the canonical symplectic 2-form. The dynamics is determined by
$$
(x, u, e, y)(t) \in \L^p \;, \quad (\dot x , e)(t) \in \D^c
$$
In local coordinates, 
\begin{gather*}
    \dot x (t) = J_{x(t)} \left( 
        \frac{\partial H}{\partial x}(x(t)) + u^i(t) \frac{\partial G_i}{\partial x}(x(t))
        \right) \\
        y_i= G_i(x(t)) \;, \quad i=1,\dots, m
\end{gather*}
The map ${J:T^*T^*Q\to TT^*Q}$ is a vector bundle isomorphism given as $J=(\Omega^{\flat})^{-1}$ where $\Omega^{\flat}:TT^*Q\to T^*T^*Q$ is defined as $X \mapsto \Omega(X, \cdot)$.
The power balance is 
$$\langle e(t)|\dot x(t) \rangle = 0 \;$$
Note that $e$ is given in terms of relation with $x$, $u$ and $y$ as dictated by $\L^p$. To inspect the relation more closely, consider the power balance in local coordinates, which boils down to
\begin{align*}
\frac{\partial H}{\partial x^i}(x(t)) \dot x^i (t) + u^i(t) \frac{\partial G_i}{\partial x}(x(t)) \dot x (t) =0    
\end{align*}
Using the chain rule, we have
$$
\frac{\d H(x(t))}{\d t} + u^i(t)\frac{\d y^i(x(t))}{\d t} = 0
$$
where $y_i=G_i$. The conjugated pair $(u,y)$ is an open energy port. The Dirac structure is closed, we do not have power ports, but there is an additional duality pairing appearing in the power balance.

\subsection{Constrained Hamiltonian systems}

Systems that are subject to constraints are confined to a part of the state space. While these constraints might arise in different ways, we look more closely at systems with gauge degrees of freedom. That is, systems with more degrees of freedom then equations that determine them. We show that Dirac's extended Hamiltonian~(\cite{dirac2001lectures}) is a Morse family linear in parameters 
provided that $0$ is a regular point of all constraints.  

Gauge theory is characterized by a non-invertible partial Legendre transformation of the Lagrangian of the system.
Let $Q$ be a smooth $d$-dimensional manifold and let the Lagrangian
${L:TQ\to\R}$ be a smooth function. The partial Legendre transform ${l:TQ\to T^*Q}$ is defined as
$$
l(q,v)=(q, \hat{p}\coloneqq \partial_vL(q,v))
$$
where $(q,v)\in \R^d\times\R^d$ denotes coordinates on the tangent bundle $TQ$.
The Hamiltonian $H_c:T^*Q\to\R$ is a function such that
$$
H_c \circ l = E
$$
where
$E:TQ\to \R$ is the Lagrangian energy function ${E(q,v)\coloneqq\hat{p}_a(q,v)v^a-L(q,v)}$. 
If an inverse $l^{-1}$ exists, i.e., if one can solve $p=\hat{p}(q,v)$ for $v$, so that $v=\hat{v}(q,p)$ where $\hat{v}$ is a function of $q$'s and $p$'s, then
${H_c (q,p)=(E\circ l^{-1})(q,p) = p_a\hat{v^a}(q,p)-L(q,\hat{v}(q,p))}$. According to the Inverse Function Theorem, $l$ is invertible at ${(q,v) \in TQ}$ if and only if the $d\times d$ matrix $W$,
$$
W_{ij}(q,v) \coloneqq \frac{\partial \hat{p}_a}{\partial v^b}(q,v) =\frac{\partial^2 L}{\partial v^i \partial v^j}(q,v)
$$
is non-singular, i.e., $\det W (q,v) \neq 0$. We assume that $W$ is singular and that $\mathrm{rank}(W)=d-N$ on all $TQ$. As a consequence of Constant Rank Theorem, 
there exists $N$ functions $\phi_j:T^*Q\to \R$ such that $\phi_j(q,p)=0$ for every $(q,p)\in l_{TQ}$ and every $j=1,\dots,N$. That is, the image $l_{TQ}$ is the zero level set of $\phi_1,\dots,\phi_N$ called primary constraint surface in $T^*Q$ and the functions are called primary constraints. Consequently, the Hamiltonian function is defined uniquely only on $l_{TQ}$, since for any $H':T^*Q\to\R$ and any $\la \in \R^k$ such that 

$$H'(x)\coloneqq H_c(x) + \la^j\phi_j(x)$$

it holds that $H' \circ l = H_c\circ l $, where $x=(q,p)$. We work with
 $H'$ since it explicitly expresses the ambiguity in defining a Hamiltonian function on entire state space $T^*Q$. This ambiguity is then passed on to equations of motion~(\cite{dirac2001lectures}), 
\begin{equation}
    \begin{gathered}
        \dot x(t) = J_{x(t)}\left( 
        \frac{\partial H_c}{\partial x}(x(t)) + \la^j(t) \frac{\partial \phi_j}{\partial x}(x(t))
        \right) \\
        \phi_j(x(t))=0 \;, \quad j=1,\dots, N
    \end{gathered}
\end{equation}
where $x:I\subseteq\R\to T^*Q$ is a differentiable curve, $\Omega$ is the canonical symplectic 2-form on $T^*Q$ and ${J:T^*T^*Q\to TT^*Q}$ is a vector bundle isomorphism given as $J=(\Omega^{\flat})^{-1}$ where $\Omega^{\flat}:TT^*Q\to T^*T^*Q$ is defined as $X \mapsto \Omega(X, \cdot)$. Note that algebraic and differential equations are intertwined, as the constraints are imposed at every point along the trajectory. We want to disentangle these equations in a sense that we want the above dynamics to take an initial condition satisfying the constraints to a later state that still satisfies all constraints. This can be achieved by requiring that the constraints do not evolve along the trajectory, i.e., $\dot \phi_j(x(t))=0$ for every $j=1,\dots,N$, which amounts to
$$
0 =  \{\phi_j, H_c \} \circ x(t) + \la^m \{ \phi_j, \phi_m \} \circ x(t)
$$
where $\{\cdot, \cdot \}$ are Poisson brackets on $(T^*Q, \Omega)$. The procedure that follows is called Dirac-Bergmann algorithm~(\cite{PONS2005491}). The condition above can either be identically satisfied, lead to new constraints (called secondary constraints) if the second term vanishes, or be used to determine the respective $\lambda^j$. If secondary constraints arise then the condition is checked again for the new list of constraints, now consisting of both primary and secondary constraints. The process is iterated until no new constraints arise. The Hamiltonian 
to which we add the final list of constraints and substitute the expressions for $\la$'s that we could determine is called the total Hamiltonian $H_T:T^*Q\to \R$
\begin{gather*}
    H_T(x) = H(x) + \la^i \phi_i(x) \\
    \phi_i(x)=0 \;, \quad i =1, \dots, k
\end{gather*}
The function $H$ contains all the terms in $H'$ for which we could determine the respective $\lambda^j$. The constraints in the final list are so-called first class constraints. If all $\lambda^j$'s can be determined, then $H_T=H$ and the dynamics of the system takes place on a symplectic submanifold of $T^*Q$~(\cite{dirac2001lectures}). If some $\lambda^j$'s are not determined, then the system has gauge degrees of freedom. 

Since $\la$'s are arbitrary, we henceforth view $H_T$ as a function on $T^*Q$ parameterized by $\la \in \R^k$. The critical set $\Sigma_{H_T}$ is
$$
\Sigma_{H_T} = \{ 
(x, \la) \in \R^{2d} \times \R^k \mid \phi_i(x)=0 \;, \quad i=1,\dots, k 
\}
$$
For $H_T$ to be a Morse family it must satisfy condition~(\ref{eq:rank}). By direct calculation one verifies that this is the case if at every point $x\in T^*Q$ at which $\phi_i(x)=0$ for $i=1,\dots, k$  and arbitrary $\la$,
$$
\mathrm{rank}
\begin{pmatrix}
    \frac{\partial \phi_1}{\partial x}(x) \\
    \vdots \\
    \frac{\partial \phi_k}{\partial x}(x) \\
\end{pmatrix} =k
$$
This is the condition for $0\in\R^k$ to be a regular value of $(\phi_1,\dots,\phi_k)$. The Lagrangian submanifold  $\L_{H_T}\subset T^*T^*Q$ generated by $H_T$ is
\begin{gather*}
    \L_{H_T} = \bigg\{
    (x,e)\in \R^{2n}\times \R^{2n} \mid \exists \la \in \R^k :\\
    e = \frac{\partial H}{\partial x}(x) + \la^i\frac{\partial \phi_i}{\partial x}(x)
    \;, \\
    \frac{\partial H_T}{\partial \la^i} = \phi_i(x) = 0 \;, \quad i =1,\dots, k \bigg\}
\end{gather*}
A system with gauge degrees of freedom is modelled by a port-Hamiltonian system $(T^*Q, \L_{H_T}, \D^c)$ where ${\D^c \subset TT^*Q\oplus T^*T^*Q}$ is a closed Dirac structure given as a graph of a Poisson bi-vector field generated by the Poisson bracket on $(T^*T^*Q, \omega)$ where $\omega$ is the canonical symplectic 2-form. 
The dynamics is determined by 
\begin{gather*}
    (x, e)(t) \in \L_{H_T} \;, \quad 
    (\dot x , e)(t) \in \D^c
\end{gather*}
As an example of a system with gauge degrees of freedom we consider a massive relativistic particle.
\begin{exmp}
    The Lagrangian $L:\R^4\to\R$ corresponds to the relativistic kinetic energy of a particle with mass $m$,
    $$
    L(q,v) = m\sqrt{\eta_{ab}v^av^b}
    $$
    where $\eta=\mathrm{diag}(-1,1)$ is the Minkowski metric. The partial Legendre transform $l:\R^4\to\R^{4*}$ given by ${\hat{p}_a\coloneqq\partial_{v^a}L}$, $a=\{0,1\}$ is not invertible. By direct calculation one can check that $\mathrm{rank}(W)=1$, so we have one constraint.
    $$
    p_a = -m\frac{\eta_{ab}v^b}{\sqrt{\eta_{cd}v^cv^d}}
    \implies \eta^{ab}p_ap_b = m^2
    $$
    Thus, $\phi(q,p)=\eta^{ab}p_ap_b-m^2$. The Lagrangian energy $E(q,v)=0$, so $H_c=0$. Since we only have one constraint we can immediately write down the total Hamiltonian
    $$
    H_T(q,p)=\la(\eta^{ab}p_ap_b-m^2)
    $$
\end{exmp}
which is a Morse family parameterized by $\la \in \R$.
The dynamics is described by a curve $x$ on $\R^4$ for which
\begin{gather*}
(x, e)(t) \in \L_{H_T} \implies e_a^q(t)=0\;, e_a^p(t)=2\la p_a(t) \;,\\ \eta^{ab}p_a(t)p_b(t)-m^2=0 \\
(\dot x , e)(t) \in \D^c \implies \dot p_a(t)=0 \;, \dot q(t) = 2\la p_a(t) \;, \\
\eta^{ab}p_a(t)p_b(t)-m^2=0
\end{gather*}
where $x(t)=(q(t),p(t))$.

\subsection{Constrained input-output Hamiltonian systems}
A restricted Morse family therefore describes a constrained Hamiltonian system with inputs and conjugated outputs. 
Considering the preceding subsections, an input-output Hamiltonian system with gauge degrees of freedom that is linear in inputs is described by a restricted Morse family $H_T:T^*Q\times\R^m\times\R^k$ on $T^Q\times\R^m$ if $0\in\R^k$ is a regular value for all constraints and there is no redundancy in inputs.
\begin{gather*}
    H_T(x, u, \la) = H(x) + u^jG_j(x) + \la^i\phi_i(x) \\
    y_j = G_j(x) \;, \quad j=1,\dots,m \\
    \phi_i(x)=0 \;, \quad i=1,\dots, k
\end{gather*}
The port-Hamiltonian system then consists of a Lagrangian submanifold that is generated by $H_T$ and a closed Dirac structure that is given as a graph of the Poisson bi-vector~(\cite{brockett1977control}). 

We now consider a massive relativistic particle in external electromagnetic field as an example of a constrained input-output Hamiltonian system. The external electromagnetic field is treated as an input.
\begin{exmp}
    The Lagrangian $L:\R^4 \to \R$ describes a relativistic particle of mass $m$ and charge $e$ interacting with electromagnetic field $A=(A_0, A_1)$ 
    $$
    L(q,v) = m \sqrt{\eta_{ab}v^av^b}+ eA_cv^c
    $$
    where $\eta=\mathrm{diag}(-1,1)$ is the Minkowski metric. The partial Legendre transform is not invertible and leads to one constraint.
    $$
    p_a = -m\frac{\eta_{ab}v^b}{\sqrt{\eta_{cd}v^cv^d}} + A_a
    \implies \eta^{ab}(p_a - eA_a)(p_b - eA_b)  = m^2
    $$
    Hence, $\phi(q,p) = \eta^{ab}(p_a - eA_a)(p_b - eA_b)  - m^2$. The Lagrangian energy $E=0$, so $H_c=0$. The total Hamiltonian is then $H_T:T^*Q\times\R^2\times\R\to\R$
    $$
    H_T(q,p, A, \la) = \la(\eta^{ab}(p_a - eA_a)(p_b - eA_b)  - m^2)
    $$
    $H_T$ is the restricted Morse family on $T^*Q\times\R^2$.
\end{exmp}




\bibliography{ifacconf}             
                                                   









\end{document}